\documentclass[12pt]{amsart}
 \usepackage{amscd,amsmath,amsthm,amssymb}
 \def\NZQ{\Bbb}               
 \def\NN{{\NZQ N}}
 \def\QQ{{\NZQ Q}}
 \def\ZZ{{\NZQ Z}}

 \def\FF{{\NZQ F}}

 \def\ab{{\bold a}}
 
 \def\xb{{\bold x}}

 \def\Xb{{\bold X}}
 
 \def\vb{{\bold v}}
 \def\opn#1#2{\def#1{\operatorname{#2}}} 
 \opn\chara{char} \opn\length{\ell} \opn\pd{pd} \opn\rk{rk}
 \opn\projdim{proj\,dim} \opn\injdim{inj\,dim} \opn\rank{rank}
 \opn\depth{depth} \opn\grade{grade} \opn\height{height}
 \opn\embdim{emb\,dim} \opn\codim{codim}
 
 \opn\Tr{Tr} \opn\bigrank{big\,rank}
 \opn\superheight{superheight}\opn\lcm{lcm}
 \opn\trdeg{tr\,deg}
 \opn\reg{reg} \opn\lreg{lreg} \opn\ini{in} \opn\lpd{lpd}
 \opn\size{size} \opn\sdepth{sdepth}
 \opn\link{link}\opn\fdepth{fdepth}\opn\lex{lex}
 \opn\bight{bight}
 \opn\div{div} \opn\Div{Div} \opn\cl{cl} \opn\Cl{Cl}
 \opn\Spec{Spec} \opn\Supp{Supp} \opn\supp{supp} \opn\Sing{Sing}
 \opn\Ass{Ass} \opn\Min{Min}\opn\Mon{Mon}
 \opn\Ann{Ann} \opn\Rad{Rad} \opn\Soc{Soc}
 \opn\Im{Im} \opn\Ker{Ker} \opn\Coker{Coker} \opn\Am{Am}
 \opn\Hom{Hom} \opn\Tor{Tor} \opn\Ext{Ext} \opn\End{End}
 \opn\Aut{Aut} \opn\id{id}
 
 \opn\nat{nat}
 \opn\pff{pf}
 \opn\Pf{Pf} \opn\GL{GL} \opn\SL{SL} \opn\mod{mod} \opn\ord{ord}
 \opn\Gin{Gin} \opn\Hilb{Hilb}\opn\sort{sort}
 \opn\Tot{Tot}
 \opn\conv{conv}
 \opn\aff{aff} \opn
 \con{conv} \opn\relint{relint} \opn\st{st}
 \opn\lk{lk} \opn\cn{cn} \opn\core{core} \opn\vol{vol}
 \opn\link{link} \opn\star{star}\opn\lex{lex}\opn\set{set}
 \opn\type{type}
 \opn\gr{gr}
 
 \def\pot#1#2{#1[\kern-0.28ex[#2]\kern-0.28ex]}

 \opn\dirlim{\underrightarrow{\lim}}
 \opn\inivlim{\underleftarrow{\lim}}
 \let\union=\cup
 \let\sect=\cap

 \def\Implies{\ifmmode\Longrightarrow \else
         \unskip${}\Longrightarrow{}$\ignorespaces\fi}
 \def\implies{\ifmmode\Rightarrow \else
         \unskip${}\Rightarrow{}$\ignorespaces\fi}
 \def\iff{\ifmmode\Longleftrightarrow \else
         \unskip${}\Longleftrightarrow{}$\ignorespaces\fi}

 \let\:=\colon
 \newtheorem{Theorem}{Theorem}[section]
 \newtheorem{Lemma}[Theorem]{Lemma}
 \newtheorem{Corollary}[Theorem]{Corollary}
 \newtheorem{Proposition}[Theorem]{Proposition}

 \newtheorem{Example}[Theorem]{Example}
 
 \newtheorem{Definition}[Theorem]{Definition}

 \let\epsilon\varepsilon
 \let\kappa=\varkappa
 \def\qed{\ifhmode\textqed\fi
       \ifmmode\ifinner\quad\qedsymbol\else\dispqed\fi\fi}
 \def\textqed{\unskip\nobreak\penalty50
        \hskip2em\hbox{}\nobreak\hfil\qedsymbol
        \parfillskip=0pt \finalhyphendemerits=0}
 \def\dispqed{\rlap{\qquad\qedsymbol}}
 
 \opn\dis{dis}
 \def\pnt{{\raise0.5mm\hbox{\large\bf.}}}
 
 \opn\Lex{Lex}

\begin{document}

 \title{Combinatorics of the integral closure of edge ideals of strong quasi-$n$-partite graphs}

\author {Monica La Barbiera and Roya Moghimipor}

\address{Monica La Barbiera, Department of Mathematics and Informatics, University of Messina, Viale Ferdinando
Stagno d’Alcontres 31, 98166 Messina, Italy} \email{monicalb@unime.it}

\address{Roya Moghimipor, Department of Mathematics, Safadasht Branch, Islamic
Azad University, Tehran, Iran} \email{roya\_moghimipour@yahoo.com}

\begin{abstract}
Combinatorial properties of some ideals related to strong quasi-$n$-partites graphs are examined.
We prove that the edge ideal of a strong quasi-$n$-partite graph is not integrally closed and
we give an expression for its integral closure. Moreover, we are able to determine the structure of the ideals of vertex covers for the edge ideals associated to a strong quasi-$n$-partite graph.

\end{abstract}


\maketitle

\section*{Introduction}
In the present paper we consider classes of monomial ideals that can arise from graph theory \cite{BM, GR,G, Ha}. More precisely, we consider classes of $n$-partite graphs and study combinatorial properties of them.
Let $G$ be a graph on the vertex set $V(G)$.
A graph $G$ is said to be $n$-partite graph if its vertex set $V=V_1\union V_2\union \cdots\union V_n$ and $V_i=\{x_{i1},\ldots,x_{im_i}\}$ for $i=1,\ldots,n$, every edge joins a vertex of $V_{i}$ with a vertex of $V_{i+1}$.
When $n=2$ these are the bipartite graphs, and when $n=3$ they are called the tripartite graphs.

A quasi-$n$-partite graph is an $n$-partite graph having loops, a strong quasi-$n$-partite graph is a complete $n$-partite graph having loops in all its vertices.
A great deal of knowledge on the strong quasi-bipartite graphs is accumulated in several papers \cite{MM,MMg,MMb, MMc, MMp}.

Algebraic objects attached to $G$ are the edge ideals $I(G)$. If $G$ is a bipartite graph having bipartition $\{V_{1},V_{2}\}$, edge ideals are monomial ideals of a polynomial ring in two sets of variables associated to such bipartition.

In detail, we are interested to handle the integral closure of edge ideals of strong quasi-$n$-partite graphs and to study algebraic aspects of it. The integral closure is a combinatoric object associated to the edge ideal.

Let $K$ be a field and $S=K[x_1,\dots,x_n]$ the polynomial
ring in $n$ variables over $K$ with each $x_i$ of degree $1$.

Next we consider the polynomial ring $T$ over $K$ in the variables
\[
x_{11},\ldots,x_{1m_1},x_{21},\ldots,x_{2m_2},\ldots,x_{n1},\ldots,x_{nm_n},
\]
and let $L$ be a monomial ideal of $T$.

We use the notation $\Xb$ for the set
$
\{x_{11},\ldots,x_{1m_1},x_{21},\dots,x_{2m_{2}},\ldots,x_{n1},\ldots,x_{nm_n}\}.
$
The integral closure $\overline{L}$ of $L$ is the set of all elements of $T$ which are integral over $L$.
The integral closure of a monomial ideal is again a monomial ideal,
$\overline{L}=(f\mid f \quad \text{is a monomial in T and} \quad f^{k}\in L^{k}, \text{for some} \quad k\geq 1)$.
Then $L$ is integrally closed, if $L=\overline{L}$.
Put $\ab_{j}=(a_{{j}_{11}},\dots,a_{{j}_{1m_{1}}},\dots,a_{{j}_{n1}},\dots,a_{{j}_{nm_{n}}})\in\NN^{m_{1}}\oplus \cdots \oplus \NN^{m_{n}}$.
If $L$ is generated by monomials $\Xb^{\ab_1},\ldots, \Xb^{\ab_q}$, a combinatorial description
of the integral closure of $L$ is the following:
$\overline{L}=(\Xb^{\lceil\alpha\rceil} \mid \alpha\in \conv(\ab_{1},\ldots,\ab_{q}))$,
where $\alpha \in \QQ^{m_{1}+\dots+m_{n}}_{+}$,
$\conv(\ab_{1},\dots,\ab_{q})$ is the set of all convex combinations of $\ab_{1},\dots,\ab_{q}$, and $\lceil\rceil$ denotes the upper right corner.
When $\Xb^{\ab}=x_{11}^{a_{11}} \cdots x_{1m_1}^{a_{1m_{1}}} \cdots x_{n1}^{a_{n1}} \cdots x_{nm_n}^{a_{nm_{n}}}$,
we write $\log(\Xb^{\ab})$ to indicate $\ab=(a_{11},\dots,a_{1m_{1}},\dots,a_{n1},\dots,a_{nm_{n}})\in \ZZ^{m_{1}+\dots+m_{n}}_{+}$.
Given a set $\mathcal{F}$ of monomials, the log set of $\mathcal{F}$, denoted by $\log(\mathcal{F})$, consists of all $\log(\Xb^{\ab})$ such that
$\Xb^{\ab}\in \mathcal{F}$.

The present paper is organized as follows. In Section~\ref{one} we study the integral closure of special classes of monomial ideals.
In Proposition \ref{graph} we prove that the integral closure of $I(G)$ of a strong quasi-$n$-partite graph $G$ is generated by binomial of degree 2.
Furthermore, we prove that the edge ideal of a strong quasi-$n$-partite graph $G$ is not integrally closed and
we give an expression for its integral closure, see Corollary \ref{generator} and Theorem \ref{structure}.

In Section \ref{two} algebraic and homological invariants of the integral closure of $I(G)$ of a strong quasi-$n$-partite graph $G$ are studied.
There is a one to one correspondence between minimal vertex covers of any graph and minimal prime ideals of its edge ideal, we generalize
to a complete $n$-partite graph with loops the notion of ideal of (minimal) vertex covers and determine the structure of the ideals of vertex covers $I_{c}(G)$ for the edge ideals associated to $G$.

In Corollary \ref{resolution} we show that the integral closure of $I(G)$ associated to a strong quasi-$n$-partite graph $G$ has a linear resolution.
We also give formulae for invariants of $T/\overline{I(G)}$ and $T/\overline{I_{c}(G)}$ for the classes of edge ideals associated to a strong quasi-$n$-partite graph $G$ such as dimension, projective dimension, depth, Castelnuovo-Mumford regularity.

\section{Integral closure of edge ideals}
\label{one}

Let $K$ be a field, and let $S=K[x_1,\dots,x_n]$ be the polynomial ring in $n$ variables over $K$ with each $x_i$ of degree $1$.

We use the notation $\Xb$ for the set
$
\{x_{11},\ldots,x_{1m_1},x_{21},\dots,x_{2m_{2}},\ldots,x_{n1},\ldots,x_{nm_n}\}.
$
Let $T=K[\Xb]$ be the polynomial ring over a field $K$ in the variables
\[
x_{11},\ldots,x_{1m_1},x_{21},\ldots,x_{2m_2},\ldots,x_{n1},\ldots,x_{nm_n},
\]
and let $L\subset T$ be a monomial ideal which is generated by monomials $\Xb^{\ab_1},\ldots, \Xb^{\ab_q}$.
Here
\begin{eqnarray*}
\Xb^{\ab_{j}}=
x_{11}^{a_{{j}_{11}}} \dots x_{1m_1}^{a_{{j}_{1m_{1}}}} x_{21}^{a_{{j}_{21}}} \cdots x_{2m_2}^{a_{{j}_{2m_{2}}}} \cdots x_{n1}^{a_{{j}_{n1}}} \cdots x_{nm_n}^{a_{{j}_{nm_{n}}}}
\end{eqnarray*}
for $\ab_{j}=(a_{{j}_{11}},\dots,a_{{j}_{1m_{1}}},a_{{j}_{21}}, \dots,a_{{j}_{2m_{2}}},\dots,a_{{j}_{n1}},\dots,a_{{j}_{nm_{n}}})\in\NN^{m_{1}}\oplus \cdots \oplus \NN^{m_{n}}$.

The {\em integral closure} of $L$ is the set of all elements of $T$ which are integral over $T$.
Since the integral closure of a monomial ideal is again a monomial ideal \cite[Theorem 1.4.2]{HH}, one has the following description of the integral closure of
$L$:
\[
\overline{L}=(f\mid f \quad \text{is a monomial in T and} \quad  f^{k}\in L^{k}, \text{for some} \quad k\geq 1).
\]
The ideal $L$ is {\em integrally closed}, if $L=\overline{L}$.

Let $\beta \in \mathbb{Q}_{+}^{n}$, where $\mathbb{Q}_{+}$ is the set of nonnegative rational numbers. We define the {\em upper right
corner} or {\em ceiling} of $\beta$ as the vector $\lceil \beta \rceil$ whose entries are given by ${\lceil \beta \rceil }_{i}$, where
\begin{equation*}
{\lceil \beta \rceil }_{i}=\left\{
\begin{array}{cc}
 \beta_{i}  & \text{ if } \quad \beta_{i} \in \NN \\
  \lfloor\beta_{i}\rfloor+1  & \text{ if } \quad \beta_{i} \notin\NN
\end{array}\right.
\end{equation*}
and where  $\lfloor \beta_{i} \rfloor$ stands for the integer part of $\beta_{i}$.
In addition, we denote the set $\{x_{i1},\dots,x_{im_{i}}\}$ by $\Xb_{i}$ for $i=1,\dots,n$.
Then the integral closure of $L$ is the monomial ideal:
\[
\overline{L}=(\Xb^{\lceil\alpha\rceil} \mid \alpha\in \conv(\ab_{1},\ldots,\ab_{q})),
\]
where
\[
\conv(\ab_{1},\dots,\ab_{q})=\left\{\sum _{j=1}^{q} \lambda_{j} \ab_{j} \biggm |  \sum_{j=1}^{q} \lambda_{j}=1, \lambda_{j}\in \mathbb{Q}_{+}\right\}
\]
is the set of all {\em convex combinations} of $\ab_{1},\dots,\ab_{q}$.

\begin{Example}
\label{integral}
{\em
Let $T=K[x_{11},x_{12},x_{13},x_{21},x_{22}]$ be the polynomial ring over a field $K$, and let
\[
L=(x_{11}^{2}x_{12}x_{21},x_{12}x_{13}x_{22}^{2},x_{13}x_{21}^{2})
\]
be the monomial ideal of $T$.
By using Normalize \cite{TN}, we obtain that
\begin{eqnarray*}
\overline{L}&=&(\Xb_{1}^{\alpha}\Xb_{2}^{\alpha'}\mid (\alpha,\alpha')\in \conv((2,1,0,1,0),(0,1,1,0,2),(0,0,1,2,0)))\\
&=&(\Xb_{1}^{\alpha}\Xb_{2}^{\alpha'}\mid (\alpha,\alpha')\in \{\lambda_{1}(2,1,0,1,0)+\lambda_{2}(0,1,1,0,2)+\lambda_{3}(0,0,1,2,0),\\
&&
\lambda_{1},\lambda_{2},\lambda_{3}\in \QQ_{+}, \lambda_{1}+\lambda_{2}+\lambda_{3}=1\})\\
&=&(x_{11}^{2}x_{12}x_{21},x_{12}x_{13}x_{22}^{2},x_{13}x_{21}^{2},x_{12}x_{13}x_{21}x_{22}).
\end{eqnarray*}
}
\end{Example}

The purpose of this section is to study the integral closure of monomial ideals in the polynomial ring $T=K[\Xb]$.

We recall the notions that come from the general theory for monomial ideals (see for instance \cite{V}).

\begin{Definition}
\label{nset}
Let $T=K[\Xb]$ be a polynomial ring over a field $K$ in the variables $x_{11},\ldots,x_{1m_1},\ldots,x_{n1},\ldots,x_{nm_n}$.
If
$\Xb^{\ab}=x_{11}^{a_{11}} \cdots x_{1m_1}^{a_{1m_{1}}} \cdots x_{n1}^{a_{n1}} \cdots x_{nm_n}^{a_{nm_{n}}}$,
we set
\[
\log(\Xb^{\ab})=\ab=(a_{11},\dots,a_{1m_{1}},\dots,a_{n1},\dots,a_{nm_{n}})\in \ZZ^{m_{1}+\dots+m_{n}}_{+}.
\]
Given a set $\mathcal{F}$ of monomials, the {\em log set} of $\mathcal{F}$, denoted by $\log(\mathcal{F})$, consists of all $\log(\Xb^{\ab})$, with
$\Xb^{\ab}\in \mathcal{F}$,
\[
\log(\mathcal{F})=\{\log(\Xb^{\ab})=\ab\in \ZZ^{m_{1}+\dots+m_{n}}_{+}\mid \Xb^{\ab}\in F\}.
\]
\end{Definition}

\begin{Example}
\label{log}
{\em
Let $\mathcal{F}=\{x_{11}^{2}x_{21},x_{11}^{3}x_{21}^{2},x_{12}x_{22}\}$ be a set of monomials in $T=K[x_{11},x_{12},x_{21},x_{22}]$.
The log set of $\mathcal{F}$ is
\[
\log(\mathcal{F})=\{(2,0,1,0),(3,0,2,0),(0,1,0,1)\}.
\]
}
\end{Example}

\begin{Definition}
\label{log}
Let $L$ be an ideal of $T=K[\Xb]$ generated by the set of monomials of $\mathcal{F}$. We define
\[
\log(L)=\{\log(\Xb^{\ab})= \ab\in \ZZ^{m_{1}+\dots+m_{n}}_{+}\mid \Xb^{\ab}\in L\}.
\]
\end{Definition}

\begin{Proposition}
\label{integrallog}
Let $L$ be a monomial ideal of $T=K[\Xb]$. Then
\[
\overline{L}=(\Xb^{\alpha}\mid \alpha\in \conv(\log(L))\cap \ZZ^{m_{1}+\dots+m_{n}}).
\]
\end{Proposition}

\begin{proof}
See (\cite{V}).
\end{proof}

Let $G$ be a graph on the vertex set $V(G)=\{v_{1},\dots,v_{n}\}$. We put $E(G)=\{\{v_{i},v_{j}\}\mid v_{i}\neq v_{j},  \quad v_{i},v_{j}\in V(G)\}$
the set of edges of $G$ and $\mathcal{L}(G)=\{\{v_{i},v_{i}\}\mid v_{i}\in V(G)\}$ the set of loops of $G$.
Furthermore we set $\mathcal{W}(G)=\mathcal{L}(G) \cup E(G)$.

An algebraic object attached to $G$ is the edge ideal
\[
I(G)=(x_{i}x_{j} \mid \{v_{i},v_{j}\}\in \mathcal{W}(G)),
\]
a monomial ideal of $S$.

If $\mathcal{L}(G)=\emptyset$, the graph $G$ is said {\em simple} or {\em loopless}, otherwise, if $\mathcal{L}(G)\neq\emptyset$, $G$ is a graph with loops.

A graph $G$ is said to be {\em n-partite} if its vertex set $V(G)=V_1\union V_2\union \cdots\union V_n$ and $V_i=\{x_{i1},\ldots,x_{im_i}\}$ for $i=1,\ldots,n$, every edge joins a vertex of $V_{i}$ with a vertex of $V_{i+1}$.
When $n=2$ these are the bipartite graphs.

A bipartite graph with its bipartition $V_{1}$ and $V_{2}$ is said {\em complete bipartite} if all the vertices of $V_{1}$ are joined to all the vertices of $V_{2}$.

\begin{Proposition}
\label{complete}
\cite[Corollary 12.2.9]{V}
If $G$ is a complete bipartite graph, then the edge ideal $I(G)$ is integrally closed.
\end{Proposition}

\begin{proof}
Let $x_{11},\dots,x_{1m_{1}},x_{21},\dots,x_{2m_{2}}$ be the vertex set of $G$, one may assume that the edges of $G$
are precisely the pairs of the form $\{x_{1p},x_{2q}\}$. Therefore,
\[
I(G)=I_{1,1}I_{2,1},
\]
where the ideals $I_{i,1}$ are squarefree Veronese ideals of degree one in the variables $x_{i1},\dots,x_{im_{i}}$ for $i=1,2$. Hence the result follows from \cite[Theorem 12.2.7]{V}.
\end{proof}

\begin{Definition}
\label{defineL^{*}}
Let $L=(\Xb_{1}^{a_{1_{1}}}\cdots \Xb_{n}^{a_{1_{n}}}, \Xb_{1}^{a_{2_{1}}}\cdots \Xb_{n}^{a_{2_{n}}},\cdots, \Xb_{1}^{a_{q_{1}}}\cdots \Xb_{n}^{a_{q_{n}}})$ be a monomial ideal of $T$, where
$\Xb_{1}^{a_{j_{1}}}\cdots \Xb_{n}^{a_{j_{n}}}$ stands for
\[
x_{11}^{a_{{j}_{11}}} \dots x_{1m_1}^{a_{{j}_{1m_{1}}}} x_{21}^{a_{{j}_{21}}} \cdots x_{2m_2}^{a_{{j}_{2m_{2}}}} \cdots x_{n1}^{a_{{j}_{n1}}} \cdots x_{nm_n}^{a_{{j}_{nm_{n}}}}
\]
for $j=1,\dots,q$. We define a monomial ideal $L^{*}$ of $T$ as
\[
L^{*}=\left(\prod_{i=1}^{n} \Xb_{i}^{\alpha_{i}}  \biggm | \alpha_{i}\in \conv(\log(F_{i}))\cap \ZZ^{m_{i}}\right),
\]
where $F_{i}=\{\Xb_{i}^{a_{1_{i}}},\Xb_{i}^{a_{2_{i}}},\dots,\Xb_{i}^{a_{q_{i}}}\}$ for $i=1,\dots,n$.
\end{Definition}

In general there is no inclusion relation between $\overline{L}$ and $L^{*}$.
The following result give a description of the integral closure expressed by its log set for the edge ideals of $n$-partite graphs.

\begin{Proposition}
\label{L^{*}}
Let $T=K[\Xb]$ be the polynomial ring over a field $K$, and $I(G)$ be the edge ideal associated to an n-partite $G$.
Then $I(G)=\overline{I(G)}$.

Consequently, $\overline{I(G)}\subseteq I^{*}(G)$.
\end{Proposition}

\begin{proof}
Let $G$ be be an $n$-partite graph on the vertex set $V(G)=V_1\union V_2\union \cdots\union V_n$, where $V_i=\{x_{i1},\ldots,x_{im_i}\}$ for $i=1,\ldots,n$. Furthermore, let $I(G)$ be the edge ideal of $G$ generated by the monomials $\Xb^{\ab_1},\ldots, \Xb^{\ab_p}$, where
\begin{eqnarray*}
\Xb^{\ab_{j}}=
x_{11}^{a_{{j}_{11}}} \dots x_{1m_1}^{a_{{j}_{1m_{1}}}} x_{21}^{a_{{j}_{21}}} \cdots x_{2m_2}^{a_{{j}_{2m_{2}}}} \cdots x_{n1}^{a_{{j}_{n1}}} \cdots x_{nm_n}^{a_{{j}_{nm_{n}}}}
\end{eqnarray*}
for $\ab_{j}=(a_{{j}_{11}},\dots,a_{{j}_{1m_{1}}},a_{{j}_{21}}, \dots,a_{{j}_{2m_{2}}},\dots,a_{{j}_{n1}},\dots,a_{{j}_{nm_{n}}})\in\NN^{m_{1}}\oplus \cdots \oplus \NN^{m_{n}}$.
Then by definition,
\[
\overline{I(G)}=(\{\Xb^{\lceil\alpha\rceil} \mid \alpha\in \conv(\ab_{1},\ldots,\ab_{p})\})
\]
and $\alpha=\sum _{j=1}^{p} \lambda_{j} \ab_{j}\in\QQ^{m_{1}+\dots+m_{n}}$ with $\sum_{j=1}^{p}\lambda_{j}=1$.

Since $G$ is a loopless graph on $m_{1}+\dots+m_{n}$ vertices, $I(G)=(x_{ih}x_{i'h'}\mid \{x_{ih},x_{i'h'}\}\in E(G))$ for $1\leq i\neq i'\leq n$ and
$1\leq h,h'\leq m_{i}$. Hence we set $\Xb^{\ab_{j}}=x_{ih}^{a_{j_{ih}}}x_{i'h'}^{a_{j_{i'h'}}}$, where $a_{j_{ih}}=a_{j_{i'h'}}=1$.
By using the definition of integral closure we have the following situations:

if $\lambda_{j}=1$ and $\lambda_{r}=0$ for any $1\leq j\neq r\leq p$, then
\[
\alpha=(a_{{j}_{11}},\dots,a_{{j}_{1m_{1}}},\dots,a_{{j}_{n1}},\dots,a_{{j}_{nm_{n}}})\in\NN^{m_{1}}\oplus \cdots \oplus \NN^{m_{n}}, \quad \Xb^{\lceil\alpha\rceil}=\Xb^{\ab_{j}},\quad  1\leq j\leq p;
\]

if $\lambda_{j}\in \QQ_{+}$, then
\[
\Xb^{\lceil\alpha\rceil}=\Xb^{\lfloor\alpha\rfloor+1}.
\]
Therefore the generators of $\overline{I(G)}$ as $\Xb^{\ab_{j}}$ and $\Xb^{\lfloor\alpha\rfloor+1}$, and hence $I(G)=\overline{I(G)}$.

Let $x_{ih}x_{i'h'}\in \overline{I(G)}$; since $F_{i}=\{x_{i1},\dots,x_{im_{i}}\}$, it follows from
Definition \ref{defineL^{*}} that $x_{ih}x_{i'h'}\in I^{*}(G)$.
\end{proof}

A graph $G$ with loops is said to be {\em quasi-n-partite} if its vertex set $V(G)=V_1\union V_2\union \cdots\union V_n$ and $V_i=\{x_{i1},\ldots,x_{im_i}\}$ for $i=1,\ldots,n$, every edge joins a vertex of $V_{i}$ with a vertex of $V_{i+1}$, and some vertices in $V(G)$ have loops.
A quasi-n-partite graph $G$ is called {\em strong} if it is a complete n-partite graph and all its vertices have loops.

Now we consider graphs with loops which edge ideals are not integrally closed and we compute the integral closure.

\begin{Proposition}
\label{graph}
Let $T=K[\Xb]$ be the polynomial ring over a field $K$, and $I(G)$ be the edge ideal of a strong quasi-n-partite $G$.
Then $\overline{I(G)}$ is generated by binomial of degree 2.
\end{Proposition}

\begin{proof}
Let $G$ be a strong quasi-n-partite with the vertex set $\Xb$, let $I(G)$ be its edge ideal.
Let $\Xb^{\ab_1},\ldots, \Xb^{\ab_p}$ be the generators of $I(G)$, where
\begin{eqnarray*}
\Xb^{\ab_{j}}=
x_{11}^{a_{{j}_{11}}} \dots x_{1m_1}^{a_{{j}_{1m_{1}}}} x_{21}^{a_{{j}_{21}}} \cdots x_{2m_2}^{a_{{j}_{2m_{2}}}} \cdots x_{n1}^{a_{{j}_{n1}}} \cdots x_{nm_n}^{a_{{j}_{nm_{n}}}}
\end{eqnarray*}
for $\ab_{j}=(a_{{j}_{11}},\dots,a_{{j}_{1m_{1}}},a_{{j}_{21}}, \dots,a_{{j}_{2m_{2}}},\dots,a_{{j}_{n1}},\dots,a_{{j}_{nm_{n}}})\in\NN^{m_{1}}\oplus \cdots \oplus \NN^{m_{n}}$.

Using the geometric description of the integral closure of a monomial ideal (\cite[Proposition 12.1.4]{V}), we have
\[
\overline{I(G)}=(\{\Xb^{\lceil\alpha\rceil} \mid \alpha\in \conv(\ab_{1},\ldots,\ab_{p})\}),
\]
with
\[
\conv(\ab_{1},\dots,\ab_{p})=\left\{\sum _{j=1}^{p} \lambda_{j} \ab_{j} \biggm |  \sum_{j=1}^{p} \lambda_{j}=1, \lambda_{j}\in \mathbb{Q}_{+}\right\}.
\]
Let $f$ be a generator of $\overline{I(G)}$. Then $f=\Xb^{\lceil\alpha\rceil}$ with $\alpha=\sum _{j=1}^{p} \lambda_{j} \ab_{j}\in \conv(\ab_{1},\dots,\ab_{p})$,
$\sum_{j=1}^{p}\lambda_{j}=1$, $\lambda_{j}\in \mathbb{Q}_{+}$. Therefore,
\[
\alpha=\left(\sum _{j=1}^{p} \lambda_{j}a_{j_{11}},\dots,\sum _{j=1}^{p} \lambda_{j}a_{j_{1m_{1}}},\dots,\sum _{j=1}^{p} \lambda_{j}a_{j_{n1}},\dots,\sum _{j=1}^{p} \lambda_{j}a_{j_{nm_{n}}}\right)\in \QQ^{m_{1}+\dots+m_{n}}_{+}
\]
with
\[
\ab_{j}=(a_{{j}_{11}},\dots,a_{{j}_{1m_{1}}},a_{{j}_{21}}, \dots,a_{{j}_{2m_{2}}},\dots,a_{{j}_{n1}},\dots,a_{{j}_{nm_{n}}})\in\NN^{m_{1}}\oplus \cdots \oplus \NN^{m_{n}}
\]
and $a_{j_{ih}}\in\{0,1,2\}$.
The generic element of $\alpha$, $\alpha_{ih}$, $1\leq i\leq n$, $1\leq h\leq m_{i}$,
\begin{eqnarray*}
\sum _{j=1}^{p} \lambda_{j}a_{j_{ih}}&=&\lambda_{1}a_{1_{ih}}+\lambda_{2}a_{2_{ih}}+\dots+\lambda_{r}.2+\dots+\lambda_{p}a_{p_{ih}}\\
&=&2\lambda_{r}+(1-\lambda_{r})\\
&=&\lambda_{r}+1\leq 2.
\end{eqnarray*}
\end{proof}

\begin{Corollary}
\label{generator}
Let $G$ be a strong quasi-n-partite and $I(G)\subset T$ the edge ideal associated to $G$. Then $I(G)$ is not integrally closed.
\end{Corollary}

\begin{proof}
Let $G$ be a graph on the vertex set $V(G)$, and let $I(G)$ be the edge ideal of a strong quasi-n-partite graph $G$.
Let $\Xb^{\ab_1},\ldots, \Xb^{\ab_p}$ be the generators of $I(G)$, where
\begin{eqnarray*}
\Xb^{\ab_{j}}=
x_{11}^{a_{{j}_{11}}} \dots x_{1m_1}^{a_{{j}_{1m_{1}}}} x_{21}^{a_{{j}_{21}}} \cdots x_{2m_2}^{a_{{j}_{2m_{2}}}} \cdots x_{n1}^{a_{{j}_{n1}}} \cdots x_{nm_n}^{a_{{j}_{nm_{n}}}}
\end{eqnarray*}
for $\ab_{j}=(a_{{j}_{11}},\dots,a_{{j}_{1m_{1}}},a_{{j}_{21}}, \dots,a_{{j}_{2m_{2}}},\dots,a_{{j}_{n1}},\dots,a_{{j}_{nm_{n}}})\in\NN^{m_{1}}\oplus \cdots \oplus \NN^{m_{n}}$.

Now we may assume that $a_{j_{ir}}=a_{l_{ir'}}=2$ for some $1\leq j,l\leq p$, where $r\neq r'$. Then Proposition \ref{graph} implies that $x_{ir}^{2},x_{ir'}^{2}$ are generators of $I(G)$.
We consider the obtained convex hull placing $\lambda_{j}=\frac{1}{2}$ and $\lambda_{l}=\frac{1}{2}$ and $\lambda_{t}=0$ $\forall t\neq j,l$. Therefore,
\[
\alpha=(0,\dots,\underbrace{1}_{\alpha_{ir}},\dots,\underbrace{1}_{\alpha_{ir'}},\dots,0)\in \conv(\ab_{1},\dots,\ab_{p}),
\]
hence $x_{ir}x_{ir'}$ is a generator of the integral closure of $I(G)$. Then $I(G)\neq \overline{I(G)}$, as desired.
\end{proof}

\begin{Example}
\label{notintegral}
{\em
Let $G$ be a strong quasi-3-partite graph on the vertex set $V(G)=\{x_{11},x_{12},x_{21},x_{22},x_{31},x_{32}\}$. Then the edge ideal of $G$ is the ideal:
\begin{eqnarray*}
&I(G)=&(x_{11}^{2},x_{12}^{2},x_{21}^{2},x_{22}^{2},x_{31}^{2},x_{32}^{2},
x_{11}x_{21},x_{11}x_{22},x_{12}x_{21},x_{12}x_{22},
x_{11}x_{31},x_{11}x_{32},\\
&&
x_{12}x_{31},x_{12}x_{32},
x_{21}x_{31},x_{21}x_{32},x_{22}x_{31},x_{22}x_{32})\subset K[x_{11},x_{12},x_{21},x_{22},x_{31},x_{32}].
\end{eqnarray*}
A computation with Normalize (\cite{TN}) gives
\begin{eqnarray*}
&\overline{I(G)}=&(x_{11}^{2},x_{12}^{2},x_{21}^{2},x_{22}^{2},x_{31}^{2},x_{32}^{2},
x_{11}x_{21},x_{11}x_{22},x_{12}x_{21},x_{12}x_{22},x_{11}x_{12},
x_{11}x_{31},\\
&&
x_{11}x_{32},x_{12}x_{31},x_{12}x_{32},x_{21}x_{22},
x_{21}x_{31},x_{21}x_{32},x_{22}x_{31},x_{22}x_{32},x_{31}x_{32}).
\end{eqnarray*}
Therefore, $I(G)\neq \overline{I(G)}$, and hence $I(G)$ is not integrally closed.
}
\end{Example}

The structure of the integral closure of $I(G)$ associated to a strong quasi-$n$-partite graph $G$ is given in the following result.

\begin{Theorem}
\label{structure}
Let $T=K[\Xb]$ be the polynomial ring over a field $K$ and $G$ be a strong quasi-n-partite graph. Then
\[
\overline{I(G)}=\sum_{l_{i}\geq 0, \sum_{i=1}^nl_{i}=2}I_{1l_{1}}I_{2l_{2}}\dots I_{nl_{n}},
\]
where the ideals $I_{il_{i}}=(x_{i1},\dots,x_{im_{i}})^{l_{i}}$ are the monomial ideals generated by all  monomials of degree $l_{i}$ in the variables $\Xb_{i}=\{x_{i1},\ldots,x_{im_i}\}$.
\end{Theorem}

\begin{proof}
Let $G$ be a strong quasi-n-partite graph on the vertex set $V(G)=V_1\union V_2\union \cdots\union V_n$ and $V_i=\{x_{i1},\ldots,x_{im_i}\}$ for $i=1,\ldots,n$.
Let $\Xb^{\ab_1},\ldots, \Xb^{\ab_q}$ be the generators of $I(G)$, where $\Xb^{\ab_{j}}$ is a monomial ideal of degree 2, namely $x_{ir}^{2}$ or $x_{ir}x_{i'r'}$ for all $1\leq  i\neq i'\leq n$.

By the geometric description of the integral closure of a monomial ideal in \cite[Proposition 12.1.4]{V}, we have
\[
\overline{I(G)}=(\{\Xb^{\lceil\alpha\rceil} \mid \alpha\in \conv(\ab_{1},\ldots,\ab_{q})\}),
\]
where
\[
\conv(\ab_{1},\dots,\ab_{q})=\left\{\sum _{j=1}^{q} \lambda_{j} \ab_{j} \biggm | \sum_{j=1}^{q} \lambda_{j}=1, \lambda_{j}\in \mathbb{Q}_{+}\right\}.
\]
Now let $f=\Xb^{\lceil\alpha\rceil}$ be a generator of $\overline{I(G)}$ with $\alpha=\sum _{j=1}^{q} \lambda_{j} \ab_{j}\in \conv(\ab_{1},\dots,\ab_{q})$,
$\sum_{j=1}^{q}\lambda_{j}=1$, $\lambda_{j}\in \mathbb{Q}_{+}$. It then follows that
\[
\alpha=\left(\sum _{j=1}^{q} \lambda_{j}a_{j_{11}},\dots,\sum _{j=1}^{q} \lambda_{j}a_{j_{1m_{1}}},\dots,\sum _{j=1}^{q} \lambda_{j}a_{j_{n1}},\dots,\sum _{j=1}^{q} \lambda_{j}a_{j_{nm_{n}}}\right)\in \QQ^{m_{1}+\dots+m_{n}}_{+},
\]
with $\ab_{j}=(a_{{j}_{11}},\dots,a_{{j}_{1m_{1}}},a_{{j}_{21}}, \dots,a_{{j}_{2m_{2}}},\dots,a_{{j}_{n1}},\dots,a_{{j}_{nm_{n}}})\in\NN^{m_{1}}\oplus \cdots \oplus \NN^{m_{n}}$.
By definition of $I(G)$, in each generator $\Xb^{\ab_{j}}$ we have $a_{j_{ip}}=0,1,2$.
We set
\[
\mathcal{M}[i]=m_{1}+\dots+m_{i-1}+m_{i+1}+\dots+m_{n}
\]
for every $i=1,\dots,n$.
Thus
\[
\sum _{j=1}^{q} \lambda_{j}a_{j_{ih}}=\lambda_{j_{it_{1}}}+\dots+\lambda_{j_{it_{\mathcal{M}[i]}}}+2\lambda_{j_{ih}},
\]
such that
$1\leq \{it_{1}\}< \{it_{2}\}< \cdots <\{it_{\mathcal{M}[i]}\}\leq q$ and $\{ih\}\neq \{it_{1}\}, \{it_{2}\}, \dots ,\{it_{\mathcal{M}[i]}\}$.

If $\lambda_{j}\in \NN$ with $\sum_{j=1}^{q}\lambda_{j}=1$ we obtain that $\Xb^{\lceil\alpha\rceil}=\Xb^{\ab_{j}}$, $\forall 1\leq j\leq q$,
that is $\Xb^{\lceil\alpha\rceil}=x_{ir}^{2}$ or $\Xb^{\lceil\alpha\rceil}=x_{ir}x_{i'r''}$ for $1\leq  i\neq i'\leq n$.

On the other hand, if $\lambda_{j_{ih}}=\frac{1}{2}$ with $\sum_{j}\lambda_{j}=1$, it follows that the monomials $\Xb^{\lceil\alpha\rceil}$ match
$x_{ir}x_{ir'}$, where $r< r'$ for all $1\leq r \leq m_{r}$ and $1\leq r'\leq m_{r'}$.

Otherwise, if $\lambda_{j}\in \QQ_{+}\backslash \NN$ one obtains a monomial $\Xb^{\lceil\alpha\rceil}$ with $\lceil\alpha\rceil \geq \ab_{j}$,
that is $\alpha_{ir}\geq a_{j_{ir}}$.
Hence the minimal system of generators of $\overline{I(G)}$ is $\{x_{ir}^{2},x_{ir}x_{ir'},x_{ir}x_{i'r''}\}$ for all $1\leq  i\neq i'\leq n$, where $r\neq r'$. Then the assertion follows.
\end{proof}

\begin{Definition}
\label{new}
Let $L$ be a monomial ideal of $T$ generated by the monomials
\[
\Xb_{1}^{a_{1_{1}}}\cdots \Xb_{n}^{a_{1_{n}}}, \Xb_{1}^{a_{2_{1}}}\cdots \Xb_{n}^{a_{2_{n}}},\cdots, \Xb_{1}^{a_{q_{1}}}\cdots \Xb_{n}^{a_{q_{n}}},
\]
where
$\Xb_{1}^{a_{j_{1}}}\cdots \Xb_{n}^{a_{j_{n}}}$ stands for
$x_{11}^{a_{{j}_{11}}} \dots x_{1m_1}^{a_{{j}_{1m_{1}}}} x_{21}^{a_{{j}_{21}}} \cdots x_{2m_2}^{a_{{j}_{2m_{2}}}} \cdots x_{n1}^{a_{{j}_{n1}}} \cdots x_{nm_n}^{a_{{j}_{nm_{n}}}}$ for $j=1,\dots,q$.
We define a monomial ideal $\widetilde{L}$ of $T$ as
\[
\widetilde{L}=(\{\Xb_{1}^{\lceil \alpha_{1}\rceil}\cdots \Xb_{n}^{\lceil \alpha_{n}\rceil}  \mid \alpha_{i}\in \conv(a_{1_{i}},\dots,a_{q_{i}}) \quad \text{for}\quad i=1,\dots,n\})
\]
with $\alpha_{i}\in \QQ^{m_{i}}_{+}$.
\end{Definition}

\begin{Proposition}
\label{definition}
Let $T=K[\Xb]$ be the polynomial ring over a field $K$ and $G$ be a strong quasi-n-partite graph. Then $\overline{I(G)}\subseteq \widetilde{I(G)}$.
\end{Proposition}

\begin{proof}
Let $G$ be a strong quasi-n-partite graph on the vertex set $V(G)=V_1\union V_2\union \cdots\union V_n$ and $V_i=\{x_{i1},\ldots,x_{im_i}\}$ for $i=1,\ldots,n$.
Let $I(G)$ be the edge ideal of $G$ generated by the monomials $\Xb_{1}^{a_{1_{1}}}\cdots \Xb_{n}^{a_{1_{n}}}, \Xb_{1}^{a_{2_{1}}}\cdots \Xb_{n}^{a_{2_{n}}},\ldots, \Xb_{1}^{a_{q_{1}}}\cdots \Xb_{n}^{a_{q_{n}}}$, where
$\Xb_{1}^{a_{j_{1}}}\cdots \Xb_{n}^{a_{j_{n}}}$ stands for
\[
x_{11}^{a_{{j}_{11}}} \dots x_{1m_1}^{a_{{j}_{1m_{1}}}} x_{21}^{a_{{j}_{21}}} \cdots x_{2m_2}^{a_{{j}_{2m_{2}}}} \cdots x_{n1}^{a_{{j}_{n1}}} \cdots x_{nm_n}^{a_{{j}_{nm_{n}}}}
\]
for $j=1,\dots,q$.
Then the integral closure of $I(G)$ is the ideal:
\[
\overline{I(G)}=\left(\prod_{i=1}^{n}\Xb_{i}^{\lceil \alpha_{i}\rceil} \biggm | (\alpha_{1},\dots,\alpha_{n})\in \conv((a_{1_{1}},\dots,a_{1_{n}}),\dots,(a_{q_{1}},\dots,a_{q_{n}})\right).
\]
By definition
\[
\widetilde{I(G)}=(\{\Xb_{1}^{\lceil \alpha_{1}\rceil}\cdots \Xb_{n}^{\lceil \alpha_{n}\rceil}  \mid \alpha_{i}\in \conv(a_{1_{i}},\dots,a_{q_{i}}) \quad \text{for}\quad i=1,\dots,n\})
\]
and let $f=\Xb_{1}^{\lceil \alpha_{1}\rceil}\cdots \Xb_{n}^{\lceil \alpha_{n}\rceil}$ be a generator of $\widetilde{I(G)}$.

By hypothesis $\alpha_{i}=\sum_{j=1}^{q}\lambda_{j}a_{j_{i}}\in \conv(a_{1_{i}},\dots,a_{q_{i}})$ with $\sum_{j=1}^{q}\lambda_{j}=1$, $1\leq i\leq n$.
Then $\alpha_{i}=(\sum_{j=1}^{q}\lambda_{j}a_{j_{i1}},\dots,\sum_{j=1}^{q}\lambda_{j}a_{j_{im_{i}}})$. Furthermore, we put
$
\mathcal{M}[i]=m_{1}+\dots+m_{i-1}+m_{i+1}+\dots+m_{n}
$
for every $i=1,\dots,n$.
Hence
\[
\sum_{j=1}^{q} \lambda_{j}a_{j_{ir}}=\lambda_{j_{ip_{1}}}+\dots+\lambda_{j_{ip_{\mathcal{M}[i]}}}+2\lambda_{j_{ir}},
\]
$1\leq \{ip_{1}\}< \{ip_{2}\}< \cdots <\{ip_{\mathcal{M}[i]}\}\leq q$ and $\{ir\}\neq \{ip_{1}\}, \{ip_{2}\}, \dots ,\{ip_{\mathcal{M}[i]}\}$.

If $\lambda_{j}\in \NN$ with $\sum_{j=1}^{m}\lambda_{j}=1$ we obtain that $\Xb_{i}^{\lceil \alpha_{i}\rceil}=x_{ip}$ or $\Xb_{i}^{\lceil \alpha_{i}\rceil}=
x_{ip}^{2}$ $\forall 1\leq i\leq n$ or $\Xb_{i}^{\lceil \alpha_{i}\rceil}=1$.

On the other hand, if $\lambda_{j}\in \QQ_{+}\setminus \NN$ it follows that $\Xb_{i}^{\lceil \alpha_{i}\rceil}=x_{i1} \cdots x_{im_{i}}$.
Therefore, $\widetilde{I(G)}$ is generated by all the products of monomials $\Xb_{i}^{\lceil \alpha_{i}\rceil}$ before defined, the assertion follows.
\end{proof}

\section{Regularity and projective dimension}
\label{two}
Let as before, $K$ be a field and $T=K[\Xb]$ be the polynomial ring over $K$ in the variables
\[
x_{11},\ldots,x_{1m_1},x_{21},\ldots,x_{2m_2},\ldots,x_{n1},\ldots,x_{nm_n},
\]
and let $L\subset T$ be a monomial ideal. We denote by $G(L)$ its unique minimal set of monomial generators.

In this section we study the regularity, depth, dim and projective dimension of monomial ideals corresponding to quasi-n-partite graphs with loops.

A monomial ideal $L$ is said to have {\em linear quotients} if there is an ordering $f_{1},\dots,f_{q}$ of monomials belonging to $G(L)$ with
$\deg(f_{1}) \leq\dots \leq \deg (f_{q})$ such that the colon ideal $(f_{1},\dots,f_{j-1}):(f_{j})$ is generated by a subset of $\Xb$ for each $2\leq j \leq q$.
Let $r_{j}$ denote the number of variables which is required to generate $(f_{1},\dots,f_{j-1}):(f_{j})$,
set $\mathfrak{r}(L)=\max_{2\leq j\leq q}r_{j}$. For this topic we refer the reader to \cite[Definition 6.3.45]{V}.

\begin{Proposition}
\label{linear}
Let $T=K[\Xb]$ be the polynomial ring over a field $K$ and $G$ be a strong quasi-n-partite graph. Then $\overline{I(G)}$ has linear quotients.
\end{Proposition}

\begin{proof}
Let $G$ be a  strong quasi-n-partite on the vertex set $V(G)=V_1\union V_2\union \cdots\union V_n$, where $V_i=\{x_{i1},\ldots,x_{im_i}\}$ for $i=1,\ldots,n$. Then $G$ is a complete n-partite graph and all its vertices have loops.

Let $\overline{I(G)}$ be the integral closure of $I(G)$ with a set of minimal monomial generators
$G(\overline{I(G)})=\{f_{11}, \dots,f_{1m_{1}},\dots,f_{q1},\dots,f_{qm_{q}}\}$.
We claim that $\overline{I(G)}$ has linear quotients with respect to the ordering
\begin{eqnarray}
\label{order}
f_{11}, \dots,f_{1m_{1}},\dots,f_{q1},\dots,f_{qm_{q}}
\end{eqnarray}
of $G(\overline{I(G)})$, where $f_{i1} <_{Lex} \cdots <_{Lex} f_{im_{i}}$
by the ordering
\[
x_{11}> \cdots > x_{1i_{1}}> \cdots > x_{n1}> \cdots > x_{ni_{n}}
\]
for all $i$.

Now let $u,v\in G(\overline{I(G)})$ be two monomials such that in (\ref{order}) the monomial $u$ appears before $v$.
In order to show that $\overline{I(G)}$ has linear quotients with respect to the above mentioned order,
we must show that there exists a variable $x_{ij}$ and a monomial $w\in G(\overline{I(G)})$ such that $x_{ij}|( u/\gcd(u,v))$, in (\ref{order})
the monomial $w$ comes before $v$, and $x_{ij}=w/\gcd(w,v)$.

We define a $K$-algebra homomorphism $\varphi :T\rightarrow S$ by $\varphi(x_{ij})=x_{i}$ for all $i,j$.
We suppose that $\varphi(u)=\varphi(v)$, and $x_{ij}$ is the greatest variable with respect to the given order on the variables such that
$x_{ij}| (u/\gcd(u,v))$. Let $w\in G(\overline{I(G)})$ be the monomial ideal with $\varphi(w)=\varphi(v)$ and $w=x_{ij}\gcd(v,w)$.
Since the order of monomials in (\ref{order}) are given by the lexicographical order, it then follows that $w$ comes before $v$ in (\ref{order}).

Next assume that $\varphi(u)\neq \varphi(v)$. Let $x_{ij}$ be the variable which divides $u/\gcd (u,v)$. Hence there exists
a monomial $w\in G(\overline{I(G)})$ such that $\varphi(w)$ coming before $\varphi(v)$.

Then there exists a variable $x_{i}$
such that $x_{i} | (\varphi(u)/\gcd(\varphi(u),\varphi(v))$ and $x_{i}=\varphi(w)/\gcd(\varphi(w),\varphi(v))$.
Therefore, the monomial $w$ before $v$ in (\ref{order}) and
\[
w=x_{ij}\gcd(w,v),
\]
as desired.
\end{proof}

As an immediate consequence we obtain the following important result

\begin{Corollary}
\label{resolution}
Let $T=K[\Xb]$ be the polynomial ring over a field $K$ and $G$ be a strong quasi-n-partite graph. Then $\overline{I(G)}$ has a linear resolution.
\end{Corollary}

\begin{proof}
Follows from the  general fact that the ideals generated in the same degree with linear quotients have
a linear resolution (see \cite{HH}, Proposition 8.2.1]).
\end{proof}

A {\em vertex cover} of a monomial ideal $L\subset T$ is a subset $\mathcal{C}$ of $\Xb$ such that each $u\in G(L)$ is divided by some $x_{ij}$ on $\mathcal{C}$.
The vertex cover $\mathcal{C}$ is called {\em minimal} if no proper subset of $\mathcal{C}$ is a vertex cover of $L$.

\medskip

Now we investigate algebraic and homological invariants of $T/\overline{I(G)}$.

\begin{Lemma}
\label{height}
Let $I(G)\subset T$ be the edge ideal of a strong quasi-n-partite graph $G$. Then
$
\height(\overline{I(G)})=m_{1}+\dots+m_{n}.
$
\end{Lemma}

\begin{proof}
Let $G$ be a  strong quasi-n-partite on the vertex set $V(G)$ and $I(G)$ be its edge ideal.
In fact, $P$ is a minimal prime ideal of $I(G)$ if and only if $P=(\mathcal{C})$, for some minimal vertex cover $\mathcal{C}$ of $G$ (\cite[Proposition 6.1.16]{V}).
The minimal cardinality of the vertex covers of $\overline{I(G)}$ is $\height(\overline{I(G)})=m_{1}+\dots+m_{n}$
being $\mathcal{C}=\{x_{11},\ldots,x_{1m_1},x_{21},\ldots,x_{2m_2},\ldots,x_{n1},\ldots,x_{nm_n}\}$
a minimal vertex cover of $\overline{I(G)}$ by construction.
\end{proof}

Consider the minimal graded free resolution of $M=T/L$ as an $T$-module:
\[
\FF: \quad \quad 0 \rightarrow \bigoplus_{j}T(-j)^{b_{gj}}\rightarrow \dots \bigoplus_{j}T(-j)^{b_{1j}}\rightarrow T\rightarrow T/L\rightarrow0.
\]
The {\em Castelnuovo-Mumford regularity} of M is defined as:
\[
\reg(M)=\max\{j-i \mid b_{ij}\neq 0\}.
\]
The numbers $b_{ij}=\dim Tor_{i}(K,M)_{j}$ are called the {\em graded Betti numbers} of $M$, and $b_{i}=\sum_{j} b_{ij}$ is called
the {\em $i$th Betti number} of $M$.

\begin{Theorem}
\label{algebraic}
Let $T=K[\Xb]$ be the polynomial ring over a field $K$ and $G$ be a strong quasi-n-partite graph. Then
\begin{itemize}
\item[(\.{i})] $\dim(T/\overline{I(G)})=0$;

\item[(\.{i}\.{i})] $\projdim(T/\overline{I(G)})=m_{1}+\dots+m_{n}$;

\item[(\.{i}\.{i}\.{i})] $\depth(T/\overline{I(G)})=0$;

\item[(\.{i}v)]$\reg(T/\overline{I(G)})=1$.

\end{itemize}
\end{Theorem}

\begin{proof}
(\.{i}) Let $G$ be a strong quasi-n-partite graph on the vertex set $V(G)$ and let $I(G)$ be its edge ideal.
By \cite[Corollary 7.2.5]{V} we have
$
\dim(T/\overline{I(G)})=\dim T-\height(\overline{I(G)}).
$
Hence Lemma \ref{height} implies that $\dim(T/\overline{I(G)})=0$.

(\.{i}\.{i}) The length of the minimal free resolution of $T/\overline{I(G)}$ over $T$ is equal to $\mathfrak{r}(\overline{I(G)})+1$ (\cite[Corollary 1.6]{HTa}). Then Proposition \ref{linear} yields $\projdim(T/\overline{I(G)})=m_{1}+\dots+m_{n}$.

(\.{i}\.{i}\.{i}) By Lemma \ref{height} we conclude that $\height(\overline{I(G)})=m_{1}+\dots+m_{n}$. Therefore $\overline{I(G)}$ is Cohen-Macaulay, and hence
$
\dim(T/\overline{I(G)})=\depth(T/\overline{I(G)})=0.
$

(\.{i}v) It follows from Theorem \ref{resolution} that $\overline{I(G)}$ has a linear resolution. Then
$
\reg(\overline{I(G)})=2.
$
\end{proof}

\begin{Example}
\label{invariants}
{\em
Let $T=K[x_{11},x_{12},x_{13},x_{21},x_{22},x_{23}]$ be a polynomial ring over a field $K$. Let $G$ be a strong quasi-bipartite graph on the vertex set
$V(G)=\{x_{11},x_{12},x_{13},x_{21},x_{22},x_{23}\}$.
Then
\begin{eqnarray*}
&I(G)=&(x_{11}^{2},x_{12}^{2},x_{13}^{2},x_{21}^{2},x_{22}^{2},x_{23}^{2},
x_{11}x_{21},x_{11}x_{22},x_{11}x_{23},x_{12}x_{21},x_{12}x_{22},x_{12}x_{23},\\
&&
x_{13}x_{21},x_{13}x_{22},x_{13}x_{23}).
\end{eqnarray*}
A computation with Normalize (\cite{TN}) gives
\begin{eqnarray*}
&\overline{I(G)}=&(x_{11}^{2},x_{12}^{2},x_{13}^{2},x_{21}^{2},x_{22}^{2},x_{23}^{2},x_{11}x_{12},x_{11}x_{13},
x_{11}x_{21},x_{11}x_{22},x_{11}x_{23},x_{12}x_{21},x_{12}x_{22},\\
&&
x_{12}x_{23},x_{12}x_{13},
x_{21}x_{22},x_{21}x_{23},x_{22}x_{23},
x_{13}x_{21},x_{13}x_{22},x_{13}x_{23}).
\end{eqnarray*}
There exists a one to one correspondence between the minimal vertex covers of $G$ and the minimal prime ideals of $I(G)$.
Then the irredundant primary decomposition of $\overline{I(G)}$ is
\begin{eqnarray*}
\overline{I(G)}&=&(x_{11},x_{12},x_{13},x_{21}^{2},x_{22}^{2},x_{23}^{2},
x_{21}x_{22},x_{21}x_{23},x_{22}x_{23})\\
&&
\sect (x_{11}^{2},x_{12}^{2},x_{13}^{2},x_{21},x_{22},x_{23},
x_{11}x_{12},x_{11}x_{13},x_{12}x_{13}).
\end{eqnarray*}
Hence Lemma \ref{height} implies that $\height(\overline{I(G)})=6$.

By \cite[Corollary 1.6]{HTa} it follows that the length of the minimal free resolution of $T/\overline{I(G)}$ over $T$ is equal to $\mathfrak{r}(\overline{I(G)})+1=6$ (\cite[Corollary 1.6]{HTa}). Therefore, Theorem \ref{algebraic} yields
\begin{itemize}
\item[(\.{i})] $\dim(T/\overline{I(G)})=0$;

\item[(\.{i}\.{i})] $\projdim(T/\overline{I(G)})=6$;

\item[(\.{i}\.{i}\.{i})] $\depth(T/\overline{I(G)})=0$;

\item[(\.{i}v)]$\reg(T/\overline{I(G)})=1$.

\end{itemize}
}
\end{Example}

In the following, we compute the Betti numbers of the integral closure of $I(G)$.

\begin{Theorem}
\label{Betti}
Let $T=K[\Xb]$ be the polynomial ring over a field $K$ and $G$ be a strong quasi-n-partite graph. Then for $i \geq0$, we have
\[
b_{i}(\overline{I(G)})= \displaystyle{m_{1}+\dots+m_{n}+1\choose m_{1}+\dots+m_{n}-i-1}\displaystyle{i+1\choose i}.
\]
\end{Theorem}

\begin{proof}
Let $G$ be a  strong quasi-n-partite on the vertex set $V(G)$, and let $I(G)$ be its edge ideal. Corollary \ref{resolution} implies that $\overline{I(G)}$ has a linear resolution.
Then Eagon-Northcatt complex resolving $\overline{I(G)}$ gives the Betti numbers; see for example \cite{B}. Alternatively one can use
the Herzog-K\"{u}hl formula to obtain the Betti numbers \cite[Theorem 1]{HK}.
Next by Auslander-Buchsbaum formula, one has $\projdim(\overline{I(G)})=m_{1}+\dots+m_{n}-1$. Therefore, by
Herzog-K\"{u}hl formula one obtains that
\begin{eqnarray*}
b_{i}(\overline{I(G)})&=&\frac{(m_{1}+\dots+m_{n}+1)!}{(m_{1}+\dots+m_{n}-i-1)!(2+i)!}\times \frac{(i+1)!}{i!}\\
&=&\displaystyle{m_{1}+\dots+m_{n}+1\choose m_{1}+\dots+m_{n}-i-1}\displaystyle{i+1\choose i},
\end{eqnarray*}
as desired.
\end{proof}

Let $I\subset S=K[x_{1},\dots,x_{n}]$ be a graded ideal. We consider $S/I$ as a standard graded $K$-algebra.
We have the following (see \cite{V}):

\begin{Proposition}
\label{proj}
Let $S/I$ be a Cohen-Macaulay ring then the type of $S/I$ is equal to the last Betti number in the minimal free resolution of $S/I$ as an $S$-module.
\end{Proposition}

\begin{Proposition}
\label{type}
Let $G$ be a strong quasi-n-partite graph and $I(G)\subset T$ its edge ideal. Then
\[
\type(T/\overline{I(G)})=m_{1}+\dots+m_{n}.
\]
\end{Proposition}

\begin{proof}
By Theorem \ref{algebraic}, we have $T/\overline{I(G)}$ is Cohen-Macaulay, then
$\dim(T/\overline{I(G)})=0$.
By Auslander-Buchsbaum formula we obtain that $\projdim(T/\overline{I(G)})=m_{1}+\dots+m_{n}$.
Hence Theorem \ref{Betti} together with Proposition \ref{proj} now yields
\begin{eqnarray*}
b_{m_{1}+\dots+m_{n}-1}(\overline{I(G)})&=&\displaystyle{m_{1}+\dots+m_{n}+1\choose 0}\displaystyle{m_{1}+\dots+m_{n}-1+1\choose m_{1}+\dots+m_{n}-1}\\
&=&m_{1}+\dots+m_{n}.
\end{eqnarray*}
Then the assertion follows.
\end{proof}

\medskip

Next we want to study the ideals of vertex covers for the class of edge ideals associated to quasi-n-partite graphs.

Let $I\subset S=K[x_1,\dots,x_n]$ be a monomial ideal. The ideal of {\em (minimal) covers} of $I$, denoted by $I_{c}$, is the ideal of $S$ generated
by all monomials $x_{i_{1}}\cdots x_{i_{r}}$ such that $(x_{i_{1}},\dots, x_{i_{r}})$ is an associated (minimal) prime ideal of $I$.
Let $G$ be a graph and let $I(G)$ be its edge ideal. We define $I_{c}(G)$ the ideal of {\em vertex covers} of $I(G)$.
Then
\[
I_{c}(G)=\bigg(\bigcap _{{\{v_{i},v_{j}\}\in E(G),i\neq j}}(x_{i},x_{j})\bigg)\sect \big(x_{p}\mid \{v_{p},v_{p}\}\in \mathcal{L}(G),\quad p\neq i,j\big).
\]

\begin{Proposition}
\label{complete}
Let $G$ be a graph  with loops, and let $I_{c}(G)$ be the ideal of vertex covers of $I(G)$.
Then for all $k\geq 1$ we have
\[
\overline{I_{c}(G)^{k}}=(\{\xb^{\lceil\alpha\rceil}\mid \alpha\in \conv(k\log(I_{c}(G))\}).
\]
\end{Proposition}

\begin{proof}
Let $G$ be a graph on the vertex set $V(G)=\{x_1,\dots,x_n\}$, and let $I_{c}(G)$ be the ideal of vertex covers of $I(G)$
generated by the monomials $\xb^{\vb_{1}},\dots,\xb^{\vb_{q}}$.
We put $\log(I_{c}(G))=\{\vb_1,\dots,\vb_q\}\subset \NN^{n}$ and we define the set
\[
k\log(I_{c}(G))=\{\vb_{p_{1}}+\vb_{p_{2}}+\dots+\vb_{p_{k}}\mid 1\leq p_{1}\leq \cdots \leq p_{k} \leq q\}.
\]
We assume that $k\log(I_{c}(G))=\{\mathfrak{v}_{1},\dots,\mathfrak{v}_{r}\}$, where $\mathfrak{v}_{p}=\vb_{p_{1}}+\vb_{p_{2}}+\dots+\vb_{p_{k}}$, with
$1\leq p_{1}\leq \cdots \leq p_{k} \leq q$. Then there are $r=\displaystyle{k+q-1\choose k}$ elements in $k\log(I_{c}(G))$, and hence
\[
\conv(k\log(I_{c}(G)))=\left\{\sum_{i=1}^{r}\lambda_{i}\mathfrak{v}_{i} \biggm | \sum _{i=1}^{r}\lambda_{i}=1, \lambda_{i}\in \QQ_{+}\right\}
\]
is the convex hull of $k\log(I_{c}(G))$.

Now let $\alpha=\sum_{i=1}^{r}\lambda_{i}\mathfrak{v}_{i}$ with $\mathfrak{v}_{i}\in k\log(I_{c}(G))$, $\sum_{i=1}^{r}\lambda_{i}=1$. $\lambda_{i}\in \QQ_{+}$.
We know that $\lceil\alpha\rceil\geq \alpha$, there is $\beta\in \QQ^{n}_{+}$ such that $\lceil\alpha\rceil=\alpha+\beta$.
Then there is $0\neq h\in \NN$ so that $h\beta\in \NN^{n}$ and $h\lambda_{i}\in \NN$ for all $i$. Therefore
\[
\xb^{h\lceil\alpha\rceil}=\xb^{h\beta}\xb^{h\alpha}=\xb^{h\beta}(\xb^{\mathfrak{v}_{1}})^{h\lambda_{1}}\cdots (\xb^{\mathfrak{v}_{r}})^{h\lambda_{r}}\in (I_{c}(G)^{k})^{h}\Rightarrow \xb^{\lceil\alpha\rceil}\in \overline{I_{c}(G)^{k}}.
\]
Conversely, let $\xb^{\gamma}\in \overline{I_{c}(G)^{k}}$, that is $\xb^{h\gamma}\in (I_{c}(G)^{k})^{h}$ for some $0 \neq h\in \NN$.
There are nonnegative integers $a_{1},\dots, a_{r}$ such that
\[
\xb^{h\gamma}=\xb^{\vartheta}(\xb^{\mathfrak{v}_{1}})^{a_{1}}\cdots (\xb^{\mathfrak{v}_{r}})^{a_{r}} \quad \text{and} \quad a_{1}+\dots+a_{r}=h.
\]
It then follows that $\gamma=(\vartheta/h)+\sum_{i=1}^{r}(a_{i}/h)\mathfrak{v}_{i}$. In addition, we set $\alpha=\sum_{i=1}^{r}(a_{i}/h)\mathfrak{v}_{i}$.
By dividing the entries of $\vartheta$ by $h$ we can write $\gamma=\beta+\xi+\alpha$, where $0\leq \xi_{i}< 1$ for all $i$ and $\beta\in \NN^{n}$.
Note that $\xi+\alpha\in \NN^{n}$. This implies that $\lceil\alpha\rceil=\xi+\alpha$. Therefore $\xb^{\gamma}=\xb^{\beta}\xb^{\lceil\alpha\rceil}$, where
$\alpha\in \conv(\mathfrak{v}_{1},\dots,\mathfrak{v}_{r})$ as desired.
\end{proof}

Let $L$ be a monomial ideal of $T$. The {big height} of $L$, denoted by, $\bight(L)$ is $\max\{\height(P)\mid P\in \Ass(T/L)\}$.

\begin{Lemma}
\label{bight}
Let $G$ be a strong quasi-n-partite graph on the vertex set $V(G)$. Then
\[
\bight(\overline{I_{c}(G)})=1.
\]
\end{Lemma}

\begin{proof}
Let $G$ be a strong quasi-n-partite graph on the vertex set $V(G)=V_1\union V_2\union \cdots\union V_n$ and $V_i=\{x_{i1},\ldots,x_{im_i}\}$ for $i=1,\ldots,n$. Let $I(G)$ be the edge ideal of $G$ generated by the monomials $\Xb_{1}^{a_{1_{1}}}\cdots \Xb_{n}^{a_{1_{n}}}, \Xb_{1}^{a_{2_{1}}}\cdots \Xb_{n}^{a_{2_{n}}},\ldots, \Xb_{1}^{a_{q_{1}}}\cdots \Xb_{n}^{a_{q_{n}}}$, where
$\Xb_{1}^{a_{j_{1}}}\cdots \Xb_{n}^{a_{j_{n}}}$ stands for
\[
x_{11}^{a_{{j}_{11}}} \dots x_{1m_1}^{a_{{j}_{1m_{1}}}} x_{21}^{a_{{j}_{21}}} \cdots x_{2m_2}^{a_{{j}_{2m_{2}}}} \cdots x_{n1}^{a_{{j}_{n1}}} \cdots x_{nm_n}^{a_{{j}_{nm_{n}}}}
\]
for $j=1,\dots,q$.
We assume that $\vb=(\underbrace{1,\dots,1}_{m_{1}- \text{times}},\dots,\underbrace{1,\dots,1}_{m_{n}- \text{times}}) \in\NN^{m_{1}}\oplus \cdots \oplus \NN^{m_{n}}$ be a vector.
Later by using \cite[Proposition 12.1.4]{V}) it turns out that
$
\overline{I_{c}(G)}=(\{\Xb^{\lceil\alpha\rceil} \mid \alpha\in \conv(\vb)\}),
$
with
\[
\conv(\vb)= \lambda_{j} \vb \quad \text{with} \quad \lambda_{j}\in \mathbb{Q}_{+}.
\]
Let $f$ be a generator of $\overline{I_{c}(G)}$. Then $f=\Xb^{\lceil\alpha\rceil}$ with $\alpha= (\lambda_{j},\dots,\lambda_{j})$, $\lambda_{j}=1$.
Then $\Xb^{\lceil\alpha\rceil}=\Xb^{\vb}$,
that is $\Xb^{\lceil\alpha\rceil}=x_{11}\cdots x_{1m_{1}}\cdots x_{n1}\cdots x_{nm_{n}}$. Therefore $I_{c}(G)$ is integrally closed, and $\overline{I_{c}(G)}=I_{c}(G)$.
The maximal cardinality of the vertex covers of $I(G)$ is $\bight(I_{c}(G))=1$
being $\{x_{ij}\}$ a maximal vertex cover of $I(G)$ by construction.
\end{proof}

\begin{Theorem}
\label{invariants,cover}
Let $T=K[\Xb]$ be the polynomial ring over a field $K$ and $G$ be a strong quasi-n-partite graph. Then
\begin{itemize}
\item[(\.{i})] $\dim(T/\overline{I_{c}(G)})=m_{1}+\dots+m_{n}-1$;

\item[(\.{i}\.{i})] $\projdim(T/\overline{I_{c}(G)})=1$;

\item[(\.{i}\.{i}\.{i})] $\depth(T/\overline{I_{c}(G)})=m_{1}+\dots+m_{n}-1$;

\item[(\.{i}v)]$\reg(T/\overline{I_{c}(G)})=m_{1}+\dots+m_{n}-1$.
\end{itemize}
\end{Theorem}

\begin{proof}

$(\.{i})$ Let $G$ be a strong quasi-n-partite graph on the vertex set $V(G)$, and let $I(G)$ be its edge ideal.
The minimal cardinality of the vertex covers of $I_{c}(G)$ is $\height(I_{c}(G))=1$
being $\mathcal{C}=\{x_{il}\}$ a minimal vertex cover of $I_{c}(G)$ by construction. Therefore, \cite[Corollary 7.2.5]{V} implies that
$
\dim(T/\overline{I_{c}(G)})=m_{1}+\dots+m_{n}-1.
$

$(\.{i}\.{i})$ Using Lemma \ref{bight} and \cite[Theorem 12.6.7]{HH}, together with \cite[Corollary 6.4.20]{V} now yield
$
\projdim(T/\overline{I_{c}(G)})=\bight(\overline{I_{c}(G)})=1.
$

$(\.{i}\.{i}\.{i})$ By the Auslander-Buchsbaum formula (see \cite[Theorem 3.5.13]{V}), one has the equality
$
\depth(T/\overline{I_{c}(G)})=\dim T-\projdim(T/\overline{I_{c}(G)})
=m_{1}+\dots+m_{n}-1.
$

$(\.{i}v)$ The ideal $\overline{I_{c}(G)}$ generated in degree $m_{1}+\dots+m_{n}$. Then \cite[Proposition 8.2.1]{HH} and
\cite[Theorem 12.6.2]{HH} says that $\overline{I_{c}(G)}$ has a $m_{1}+\dots+m_{n}$-linear resolution. Therefore,
$\reg(\overline{I_{c}(G)})=m_{1}+\dots+m_{n}$,
as desired.
\end{proof}

\end{document}